\documentclass[oneside,a4paper,12pt,reqno]{amsart}

\usepackage{upref,amsxtra,amssymb}

\newtheorem{theorem}{Theorem}[section]
\newtheorem{lemma}[theorem]{Lemma}
\newtheorem{prop}[theorem]{Proposition}
\newtheorem{cor}[theorem]{Corollary}

\begin{document}

\title{Isomorphism rigidity of commuting automorphisms}
\author{Siddhartha Bhattacharya}
\address{ School of Mathematics, 
Tata Institute of Fundamental Research, Mumbai 400005, India} 
\email{siddhart@math.tifr.res.in}
\subjclass[2000]{37A35, 37A15}
\keywords{Rigidity, commuting automorphisms, entropy}

 \begin{abstract}
 Let $d > 1$, and let 
 $(X,\alpha)$ and $(Y,\beta)$ be two zero-entropy ${\mathbb{Z}}^d$-actions 
 on compact abelian groups by $d$ commuting automorphisms. We show that
 if all lower rank subactions of $\alpha$ and $\beta$ have
 completely positive entropy, then 
 any measurable equivariant map from $X$
 to $Y$ is an affine map. In particular, two such actions are
 measurably
 conjugate if and only if they are algebraically conjugate. 
 \end{abstract}

  \maketitle




\section{Introduction}

It is well known that ergodic automorphisms of compact abelian 
groups are measurably isomorphic with Bernoulli shifts (cf. e.g. 
\cite{A}). In particular, entropy is the only measurable conjugacy
invariant for such automorphisms.
On the other hand, for $d > 1$, 
mixing zero-entropy
${\mathbb{Z}}^d$-actions on compact abelian groups by $d$ commuting
automorphisms tend to exhibit remarkable rigidity properties. 
In this paper we study rigidity of measurable structure of these
actions.

Before stating our main result, we recall a few basic definitions. 
Throughout this paper the term 
{\it compact abelian group\/} will denote an infinite 
compact metrizable abelian group. 
An {\it algebraic ${\mathbb{Z}}^d$-action \/} $(X,\alpha) $ 
is an action $\alpha$ of ${\mathbb{Z}}^d$ on a 
compact abelian group $X$ by continuous
automorphisms. Any such action preserves $\lambda_{X}$,
the normalized Haar measure on $X$.
A {\it lower rank
subaction \/} of $\alpha$ is the restriction of $\alpha$ to
a subgroup $\Lambda\subset {\mathbb{Z}}^d$ with $\mbox{rank}(\Lambda) < d$.
If  $(X,\alpha )$ and $(Y,\beta) $ are 
algebraic ${\mathbb{Z}}^d $-actions then
 a Borel map $\phi \colon X\rightarrow Y$ is said to be 
{\it equivariant\/} if
$\phi\circ \alpha({\bf n}) = \beta({\bf n})\circ \phi\ $
 $\lambda_X$-a.e., for every ${\bf n}\in {{\mathbb{Z}}^d}$.
The actions $\alpha $ and $\beta $ are {\it measurably \/} or 
{\it algebraically conjugate \/} if the map $\phi $ 
 can be chosen to be a Borel isomorphism
 or a continuous group isomorphism.
A map $\psi \colon X\rightarrow Y$ is {\it affine \/} if there
exist a continuous group homomorphism 
$\psi^{'}\colon X\rightarrow Y$ 
and an element $y\in Y$ such that
	$
\psi (x)=\psi^{'}(x)+y
	$
almost everywhere with respect to $\lambda_{X}$. 

In this paper we prove the following theorem :
\begin{theorem}
\label{main}
Let $d > 1$, and let
 $(X,\alpha)$ and $(Y,\beta)$ be  
zero-entropy algebraic ${\mathbb{Z}}^d$-actions such that 
all lower rank subactions
of $\alpha$ and $\beta$ 
have  completely positive entropy. Then every measurable 
equivariant map $f : X\rightarrow Y$ is an affine map.
In particular, two such actions are 
measurably
conjugate if and only if they are algebraically conjugate. 
\end{theorem}

\medskip
Several results in this direction were obtained in recent years.
The case when both $X$ and $Y$ are zero-dimensional was studied
in \cite{BS}, \cite{E} and \cite{KS}. For a certain class of actions
on connected groups with the property that entropy
of every individual element in  ${\mathbb{Z}}^d$ is finite, in \cite{EL}
and \cite{KKS} the above
result was proved as a consequence of more general results on 
invariant measures.
 In \cite{EW}, Theorem 1.1 was proved for
a class of ${\mathbb{Z}}^3$-actions on connected groups, where 
every individual element has infinite entropy.

The more general situation where $\alpha$ and $\beta$ are
arbitrary mixing zero-entropy actions, was studied in \cite{B1},
\cite{B2} and \cite{BW}. Although these actions also exhibit several 
 rigidity properties,  there are mixing zero-entropy algebraic
${\mathbb{Z}}^d$-actions on zero-dimensional groups for which the analogue of
Theorem 1.1 fails (cf. \cite{B1, BS}).
It is not known whether one
can construct similar examples when both $X$ and $Y$ are connected. 



\section{Background}
Let
${R_d} = {{\mathbb{Z}}}[u_{1}^{\pm 1},\ldots ,u_{d}^{\pm 1}]$ be 
the ring of Laurent polynomials with integral 
coefficients in the commuting variables 
$u_1,\ldots, u_d$. We write $f\in {R_d}$ as
$$f =\sum_{\mathbf{n}\in {\mathbb{Z}}^d}f_{\mathbf{n}}u^{\mathbf{n}}$$
with $u^{\mathbf{n}}= u_1^{n_1}\cdots u_d^{n_d}$ and
$f_{\mathbf{n}}\in {{\mathbb{Z}}}$ for all ${\mathbf{n}} = (n_1,\dots
,n_d)\in {\mathbb{Z}}^d$, 
where $f_{\mathbf{n}} = 0$ for all but finitely many 
${\bf n}\in {\mathbb{Z}}^d$.

If $\alpha $ is an algebraic ${{\mathbb{Z}}^d}$-action on a compact 
abelian group $X$, then the additively-written dual group 
$M=\widehat{X}$ is a module over the ring ${R_d}$ with respect
to the operation
$$ f\cdot a = 
\sum_{{\mathbf{n}}
\in{\mathbb{Z}}^d} f_{\bf n}\ {\widehat{\alpha}}({\bf n})(a)$$
for $f\in {R_d}$ and $a\in M$, where ${\widehat{\alpha}}({\bf n})$ 
denotes the automorphism of $\widehat{X}$ dual to
$\alpha ({\bf n})$. 
The module $M=\widehat{X}$ is called the {\it dual module\/} 
of $\alpha $.

Conversely, if $M$ is a module over ${R_d}$, then we obtain an 
algebraic 
${\mathbb{Z}}^d$-action $\alpha_M$ on $X_M=\widehat{M}$ by setting
${\widehat{\alpha _M}}({\bf n})(a) = u^{\bf n}\cdot a $
for every ${\bf n}\in{\mathbb{Z}}^d$ and $a\in M$. 
Clearly, $M$ is the dual module of $\alpha _M$.

If $(X,\alpha)$ is an algebraic ${\mathbb{Z}}^d$-action and 
$\Lambda\subset {\mathbb{Z}}^d$ is a subgroup, then $\alpha^{\Lambda}$ will
denote the restriction of $\alpha$ to $\Lambda$ . 
For a ${R_d}$-module $M$, by $F(M)$ we denote the 
submodule consisting of all $m\in M$ such that
${R_d}\cdot m$ is finitely generated as an additive group.
If $N\subset M$ is a submodule then $N^{\perp}\subset X_M$
will denote the subgroup consisting of all $x$ such that
$\chi(x) = 1$ for all $\chi\in N$. By duality theory, the
correspondence $N\mapsto N^{\perp}$ is a order-reversing
bijection from the set of all submodules of $M$ to the  set of all
$\alpha_M$-invariant closed subgroups of $X_M$.

For an algebraic ${\mathbb{Z}}^d$-action $(X,\alpha)$, the topological entropy
$h_{\rm top}(\alpha)$ coincides with the metric entropy
$h_{\lambda_X}(\alpha)$. 
If an algebraic ${\mathbb{Z}}^d$-action has completely positive entropy then 
it is Bernoulli (cf. \cite{RS}).
 Recall that a prime ideal ${\mathfrak{p}}\subset {R_d}$ is
  {\it associated with \/} a ${R_d}$-module $M$ if
  ${\mathfrak{p}} = 
\mbox{ann}(a) = \{f\in {R_d} :  f\cdot a=0_M\}$ for
  some $a\in M$. The set of prime ideals associated with $M$ will 
be denoted by $\mbox{Asc}(M)$. We recall the following results
from \cite{LSW} and \cite{S}.

\begin{lemma}
\label{mix}
Let $d\ge 1$, and let $M$ be a countable ${R_d}$-module.
Then $\alpha_M$ is mixing if and only if 
$\{u^{\bf n} -1 :\ {\bf n}\in {{\mathbb{Z}}^d}\}\cap {\mathfrak{p}} = \{0\}$
for every ${\mathfrak{p}}\in\mbox{Asc}(M)$.
\end{lemma}

\begin{lemma}
\label{ent}
Let $d\ge 1$, and 
let $M$ be a countable ${R_d}$-module such that the action 
$\alpha_M$ is mixing.
\begin{enumerate}
\item{$\alpha_M$ has zero entropy if and only if every 
${\mathfrak{p}}\in\mbox{Asc}(M)$ is non-principal.   }
\item{$\alpha_M$ has completely positive entropy if and only if every 
${\mathfrak{p}}\in\mbox{Asc}(M)$ is principal.}
\item{If $M$ is Noetherian then $\alpha_M$ has
finite entropy if and only if every ${\mathfrak{p}}\in{\mbox{Asc}}(M)$
is non-zero.}
\end{enumerate}
\end{lemma}
The {\it Krull dimension \/} $\mbox{kdim}(R)$ of a ring $R$
is the length $l$ of the longest chain
$$\{ 0\} = {\mathfrak{p}}_{0}\subset {\mathfrak{p}}_{1}\subset 
\cdots \subset {\mathfrak{p}}_{l}$$
of distinct prime ideals in $R$ (see \cite[Chapter 8]{EB} for necessary 
background). The ring ${R_d}$ has Krull
dimension $d + 1$.
If ${\mathbb{K}}$ is a field then the {\it transcendence degree \/}
$\mbox{tdeg}({\mathbb{K}})$ of ${\mathbb{K}}$ is the 
maximum number of elements 
in ${\mathbb{K}}$ that are algebraically independent over the 
prime subfield
of ${\mathbb{K}}$.  If ${\mathfrak{p}}\subset{R_d}$ is a prime ideal, and
${\mathbb{K}}$ is the field of fractions of ${R_d}/{\mathfrak{p}}$, 
then $\mbox{kdim}({R_d}/{\mathfrak{p}}) = \mbox{tdeg}({\mathbb{K}}) + 1$ if
$\mbox{char}({\mathbb{K}}) = 0$, and 
$\mbox{kdim}({R_d}/{\mathfrak{p}}) = \mbox{tdeg}({\mathbb{K}})$ if 
$\mbox{char}({\mathbb{K}}) > 0$.

For a compact abelian group $X$, we denote the connected component
containing identity by $X^{0}$. A connected compact abelian group
$X$ is {\it finite-dimensional\/} if the dual group $\widehat{X}$
is isomorphic with a subgroup of ${\mathbb{Q}}^n$ for some $n\ge 1$. 
\begin{prop}
\label{dim}
Let $M$ be a  Noetherian
${R}_{2}$-module such that the action $\alpha_M$ is mixing and has
zero entropy. Then  $X^{0}_{M}$ is a finite-dimensional group.
\end{prop}
{\it Proof.\/}
Since the restriction of $\alpha_M$ to $X^{0}_{M}$ is also mixing
and has zero entropy (cf. \cite[Theorem 3.6]{S}), it is enough to consider
the case when $X_M$ is connected. Clearly, if there exists 
a ${R}_{2}$-module $N\supset M$ such that $X_N$ is 
finite-dimensional, then $X_M$ is 
finite-dimensional. Similarly, if there exists a submodule 
$M_1\subset M$ such that both $X_{M_1}$ and $X_{M/M_1}$ are
finite-dimensional then  $X_M$ is 
finite-dimensional. As $M$ is a Noetherian module, there  
exists
a  Noetherian ${R}_{2}$-module $N\supset M$ and a finite sequence of 
submodules
$$ \{ 0\} = N_0\subset N_{1}\subset\cdots \subset N_k = N,$$
such that for each $i\ge 1$,
$N_{i+1}/N_i = {{R}_{2}}/{\mathfrak{p}}$ 
for some ${\mathfrak{p}}\in\mbox{Asc}(M)$
(cf. \cite[Corollary 6.3]{S}).
Hence  we may
assume that  $M = {{R}_{2}}/{\mathfrak{p}}$ for some prime ideal 
${\mathfrak{p}}\subset {{R}_{2}}$.
Let ${\mathbb{K}}$ denote the field of fractions of ${{R}_{2}}/{\mathfrak{p}}$.
As ${\mathfrak{p}}$ is non-principal and $\mbox{char}({\mathbb{K}}) = 0$, 
$\mbox{tdeg}({\mathbb{K}}) = \mbox{kdim}({{R}_{2}}/{\mathfrak{p}}) -1 = 0$.
Since ${{R}_{2}}/{\mathfrak{p}}$ is a finitely generated ring,
this shows that as an additive group ${\mathbb{K}}$ is isomorphic 
with ${\mathbb{Q}}^n$ 
for some $n \ge 1$, which proves the given assertion.
$\hfill \Box$

\medskip
Our next lemma is a special case of Theorem \ref{main}. It 
follows from known results on invariant measures of
algebraic ${\mathbb{Z}}^d$-actions, and a joining argument used in 
\cite{KKS} and \cite{KS}. We sketch a proof for convenience of the
reader.
\begin{lemma}
\label{toral}
Let $d > 1$, and let $M$ and $N$ be countable ${R_d}$-modules with the
following properties :
\begin{enumerate}
\item{ $F(M) = M$,}
\item{The actions $\alpha_M$ and $\alpha_N$ have zero entropy,
and all lower rank subactions of $\alpha_M$ and $\alpha_N$
have completely positive entropy.}
\end{enumerate}
Then every measurable equivariant map $f$ from
$(X_M, \alpha_M)$ to $(X_N,\alpha_N)$ is an affine map.
\end{lemma}

{\it Proof.\/} We choose an arbitrary $m\in M$. Since ${R_d}\cdot m$
is finitely generated as an additive group and the action
$\alpha_M$ is mixing, as an additive group ${R_d}\cdot m$ is isomorphic
with ${\mathbb{Z}}^k$ for some 
$k\ge 0$. Hence $M$ is torsion free, and by duality
$X_M$ is connected. In particular, the action $\alpha_M$
is mixing of all orders (cf. \cite{SW}). We also observe that
$X_{{R_d}\cdot m}$ is isomorphic with a torus and
$(X_{{R_d}\cdot m},\alpha_{{R_d}\cdot m})$ is a factor of 
$(X_M, \alpha_M)$. As every lower rank
subaction of $\alpha_M$ has completely positive entropy, we deduce
that either $X_M$ is trivial or $d= 2$.

Define $i : X_M\rightarrow X_M\times X_N$ and a measure $\mu$ on 
$X_M\times X_N$ by $i(x) = (x, f(x))$, $\mu = i_{*}(\lambda_{X_M})$.
It is easy to see that $f$ is an affine map if and only if $\mu$ is
affine, i.e. $\mu$ is a translate of the Haar measure on a closed
subgroup of $X_M\times X_N$. An elementary Harmonic analysis argument
shows that $\mu$ is of this form if and only if
$|{\widehat{\mu}}(\chi)| = 0$ or $1$ for every $\chi\in M\times N$,
where ${\widehat{\mu}}$ is the Fourier transform of $\mu$.

We fix $\chi\in M\times N$ and set $P = {R_d}\cdot\chi$. Let $\pi$
denote the projection map from $X_M\times X_N$ to $X_M$, $\pi^{'}$
denote the projection map from $X_M\times X_N$ to $X_P$, and $\mu_{P}$
denote the measure $\pi^{'}_{*}(\mu)$. Since $\pi$ is a measurable
conjugacy from $(X_M\times X_N, \alpha_M\times \alpha_N, \mu)$
to $(X_M, \alpha_M, \lambda_{X_M})$ and $\pi^{'}$
is a measurable factor map from 
$(X_M\times X_N, \alpha_M\times \alpha_N, \mu)$
 to $(X_P, \alpha_P, \mu_{P})$,  with respect to the
measure $\mu_{P}$ the action $\alpha_{P}$ is mixing of all orders, and
every lower rank subaction of $\alpha_{P}$ has completely positive
entropy. As
$X_P/X_{P}^{0}$ is zero-dimensional, from \cite{SI} we deduce that
the image of $\mu_{P}$ on $X_P/X_{P}^{0}$ is concentrated at a
point. Hence $\mu_{P}$ is a translate of some  $\alpha_P$-invariant measure
$\nu$ on $X_{P}^{0}$. Since we only need to consider the case when
$d=2$, and $\alpha_{P}$ has zero-entropy with respect to $\lambda_{X_P}$,
by the previous proposition we may assume that
$X_{P}^{0}$ is a finite-dimensional group. As all lower
rank subactions of $\alpha_{P}$ has completely positive entropy
with respect to the measure $\nu$, 
 from \cite[Theorem 1.3]{EL}
we deduce that $\nu$ is an affine measure.
This implies that
$\mu_{P}$ is an affine measure and 
$|{\widehat{\mu}}(\chi)| = |{\widehat{\mu_{P}}}(\chi)| \in \{0,1\}$.
$\hfill \Box$

\section{Continuity and finiteness of entropy}

\bigskip
In this section we prove two lemmas which will be used in the proof of
Theorem \ref{main}. We begin with a few notations.
Let $(Y, d_Y)$ be a compact metric space, and let $\beta$ be an action
of ${\mathbb{Z}}^d$ on $Y$ by homeomorphisms. We denote the set of all
$\beta$-invariant probability measures on $Y$ by $M_\beta(Y)$. With
respect to the $\rm{weak}^*$ topology $M_\beta(Y)$ is a compact 
convex set. For any $\mu\in M_{\beta}(Y)$, $h_{\mu}(\beta)$ will
denote the
entropy of the action $\beta$ with respect to the measure $\mu$.
For $n\ge 1$, let $B_n\subset {\mathbb{Z}}^d$  denote the rectangle
$\{0,\ldots , n\}^d$, and let $d_n$ denote the metric on $Y$ defined by
$$ d_n(y_1, y_2) = \max_{{\bf m}\in B_n}
\ d_Y (\beta({\bf  m})(y_1), \beta({\bf  m})(y_2)).$$
For a closed set $C\subset Y$ and $\epsilon > 0$, let
$S_n(C,\epsilon)$ be the largest cardinality of an
$\epsilon$-separating
set in $C$ with respect to the metric $d_n$.
We set
$$ S(C,\epsilon) = { \limsup_{n\mapsto\infty}}{\frac{1}{n}} 
\mbox{log}
S_n(C,\epsilon),\ \  h_{{\rm top}}(\beta, C) = \lim_{\epsilon\mapsto
  0} S(C,\epsilon).$$
 Note that $ h_{{\rm top}}(\beta, Y)$ is the topological entropy of the 
action $\beta$.
 For any $y\in Y$ and $t > 0$, let $A_{t}(y)$ denote the
set of all $y^{'}\in Y$ with the property that $d_n(y,y^{'})\le t$
for all $n\ge 1$. 
The action $\beta$ is said to be 
{\it asymptotically $h$-expansive \/} if 
$\sup_{y\in Y} h_{{\rm top}}(\beta, A_{t}(y)) \mapsto 0$ as $t\mapsto 0$.

We now recall a result from \cite{M} which is a generalization
of the well known fact that for an expansive ${\mathbb{Z}}^d$-action
$\beta$, the map $\mu\mapsto h_{\mu}(\beta)$ is
upper semicontinuous (see \cite[Corollary 2.1 and Theorem 4.2]{M}).
Although in \cite{M} the results are stated for
${\mathbb{Z}}$-actions, 
the proofs 
can easily be extended to ${\mathbb{Z}}^d$-actions for any $d \ge 1$.
\begin{lemma}
\label{semi}
Let $Y$ be a compact metric space, and let $\beta$ be 
an asymptotically $h$-expansive ${\mathbb{Z}}^d$-action on $Y$. 
Then the map $\mu\mapsto h_{\mu}(\beta)$ is upper semicontinuous.
\end{lemma}
We note a simple consequence of the above result.
\begin{lemma}
\label{semia}
Let $(Y,\beta)$ be an  algebraic ${\mathbb{Z}}^d$-action 
with finite topological entropy. Then $\mu\mapsto h_{\mu}(\beta)$
is an upper semicontinuous map from $M_\beta(Y)$ to ${\mathbb{R}}$.
\end{lemma}

{\it Proof.\/} In view of the previous lemma it is enough to show that 
the action $(Y, \beta)$ is asymptotically $h$-expansive. Let $d_Y$ be
a translation invariant metric on $Y$. For $t > 0$ the closed set 
$Y_t = A_{t}(0)$ is invariant under $\beta$. Note that 
for any $y\in Y$ and $n\ge 1$, the map $x\mapsto xy$ is an isometry
from $(Y_t, d_n)$ to $(A_{t}(y), d_n)$. 
This implies that for any $y\in Y$ and $t > 0$,
$ h_{{\rm top}}(\beta, A_{t}(y))) =  h_{{\rm top}}(\beta, Y_t)$.
Suppose that  
$\mbox{lim sup}_{t\mapsto 0} h_{{\rm top}}(\beta, Y_t) = c > 0$. Since 
$h_{{\rm top}}(\beta)$ is finite, we can choose $\epsilon > 0$ such that 
$S(Y, \epsilon ) \ge h_{{\rm top}}(\beta) - c/4$. We choose 
$0 < t,\delta < \epsilon/3$ such that $h_{{\rm top}}(\beta, Y_t) > c/2$ and 
 $S(Y_t, \delta ) > c/3$.
If $Y_1$ is a $(n,\epsilon)$-separating set in $Y$
and $Y_2$ is a $(n,\delta)$ -separating set in $Y_t$, then
$Y_1\cdot Y_2$ is a $(n,\delta)$ -separating set in $Y$.
Hence 
$$S_n(Y,\epsilon)\cdot S_n(Y_{t},\delta)
\le S_{n}(Y,\delta),$$
which implies that
 $h_{{\rm top}}(\beta) - c/4 + c/3 \le h_{{\rm top}}(\beta)$. This
contradiction shows that the action $\beta$ is 
 asymptotically $h$-expansive.
$\hfill \Box$

\bigskip
We now prove a lemma on finiteness of entropy of rank $d-1$ subactions.
It is easy to construct examples of zero-entropy algebraic
${\mathbb{Z}}^d$-actions with Noetherian dual module,
which admit rank $d-1$ subactions with infinite entropy. For example,
if $A$ is the shift automorphism on ${{\mathbb{T}}}^{{\mathbb{Z}}}$ 
and $B = \mbox{Id}$,
then the ${\mathbb{Z}}^2$-action generated by $A$ and $B$ has zero entropy,
but the cyclic action generated by $A$ has infinite entropy. We show
that this situation can not occur if the action  is assumed to be mixing.

\begin{lemma}
\label{finite}
Let $d > 1$, and let 
$M$ be a  Noetherian
${R_d}$-module such that the action $\alpha_M$ is mixing and has
zero entropy. Then for any subgroup $\Lambda\subset{{\mathbb{Z}}^d}$
with $\mbox{rank}(\Lambda) = d-1$, the action
$\alpha_{M}^{\Lambda}$ has finite entropy.
\end{lemma}

{\it Proof.\/}
Let ${C}_{d}$ denote the class of  Noetherian 
${R_d}$-modules $M$ with the property that
$\alpha_M$ is a mixing action with zero-entropy, and all rank $d-1$
subactions of $\alpha_M$ have finite entropy. Since any factor of 
$\alpha_M$ also has these properties, the class ${C}_{d}$
is closed under taking submodules. If $M\subset N$ is a submodule such
that $M$ and $N/M$ lie in ${C}_{d}$, then from 
entropy addition formula (cf. \cite[Theorem 14.1]{S}) it follows that
$N\in {C}_{d}$. Now as in the proof of Proposition \ref{dim},
we may assume that  $M = {R_d}/{\mathfrak{p}}$ for some prime 
ideal ${\mathfrak{p}}\subset {R_d}$.

Let $\Lambda_0\subset {\mathbb{Z}}^d$ denote the subgroup consisting of all
${\bf n} = (n_1,\ldots ,n_d)$ with $n_d = 0$. 
We choose an injective endomorphism 
$\phi : {{\mathbb{Z}}^d}\rightarrow {{\mathbb{Z}}^d}$
such that $\phi(\Lambda_0) = \Lambda$, 
and define a ${\mathbb{Z}}^{d}$-action $\alpha^{0}$
on $X_M$ by setting $\alpha^{0}({\bf n}) = \alpha(\phi({\bf n}))$ for all
${\bf n}\in {{\mathbb{Z}}^d}$. Clearly, the action $\alpha^{0}$ is mixing and
has zero entropy. Replacing $M$ by the dual module of the action
$\alpha^{0}$ if necessary, we may assume that $\Lambda = \Lambda_0$.

Let $R_{d-1}\subset{R_d}$ denote the subring 
${\mathbb{Z}}[ u^{\pm 1},\ldots ,u^{\pm d-1}]$, and let $R$ denote the image
of $R_{d-1}$ in ${R_d}/{\mathfrak{p}} $.
Then  as a $R_{d-1}$-module $R$  is
isomorphic with ${R_{d-1}}/{R_{d-1}}\cap {\mathfrak{p}}$.
Let ${\mathbb{F}}$ denote the field of fractions of $R$, and let ${\mathbb{K}}$
denote the field of fractions of ${R_d}/{\mathfrak{p}}$. 

We claim that the
$\Lambda$-action $\alpha_{R}$ has finite entropy.
By Lemma \ref{ent} it is enough to show that 
${R_{d-1}}\cap {\mathfrak{p}} \ne \{ 0\}$.
This is obvious if $\mbox{char}({\mathbb{K}}) > 0$.
If $\mbox{char}({\mathbb{K}}) = 0$ then 
$\mbox{tdeg}({\mathbb{K}}) = \mbox{kdim}({R_d}/{\mathfrak{p}}) - 1 
\le d-2$, as 
${\mathfrak{p}}$ is a non-principal prime ideal. Since the 
transcendence degree of the field of fractions of $R_{d-1}$
is $d-1$, the map from $R_{d-1}$ to
${R_d}/{\mathfrak{p}}\subset {\mathbb{K}}$ induced by the inclusion map
$i : R_{d-1}\rightarrow {R_d}$ is not injective. This shows that 
$R_{d-1}\cap {\mathfrak{p}} \ne \{0\}$, which proves the claim.
Note that for any finitely generated 
${R_{d-1}}$-module $N\subset {\mathbb{F}}$ we can 
choose $r_{0}\in R$ such that $r_{0}\cdot N\subset R$. Hence
$N$ is isomorphic with a submodule of ${\mathbb{F}}$, i.e. $\alpha_N$ is an
algebraic factor of $\alpha_{{\mathbb{F}}}$. In particular,
$h_{{\rm top}}(\alpha_N)\le h_{{\rm top}}(\alpha_R)$.
We choose an increasing  sequence of  ${R_{d-1}}$-submodules 
$N_1\subset N_2\subset\cdots\subset {\mathbb{F}}$ such that $\cup_j N_j = 
{\mathbb{F}}$.
Since $N_{j}^{\perp}\mapsto 0_{X_{{\mathbb{F}}}}$ as $j\mapsto\infty$,
from \cite[Lemma 13.6]{S} it follows that
$$h_{{\rm top}}(\alpha_{{\mathbb{F}}}) = 
\lim_{j\mapsto \infty} h_{{\rm top}}(\alpha_{N_j})
\le h_{{\rm top}}(\alpha_{R}) < \infty.$$

Now we consider the cases $[{\mathbb{K}} : {\mathbb{F}} ] < \infty$
and 
$[{\mathbb{K}} : {\mathbb{F}} ] =
\infty$, separately. In the former case,
the
 ${R_{d-1}}$-module 
${\mathbb{K}}$ is a direct sum
of finitely many copies of ${\mathbb{F}}$. 
Since $h_{{\rm top}}(\alpha_{{\mathbb{F}}}) < \infty$ and $M = 
{R_d}/{\mathfrak{p}}$
is a $R_{d-1}$-submodule of ${\mathbb{K}}$, the action
$\alpha_{M}^{\Lambda}$ has finite entropy.
In the latter case the element $u_d$ 
is not algebraic over ${\mathbb{F}}$. This
can happen only if ${\mathfrak{p}}$ is contained in ${R_{d-1}}$. It is 
easy to see
that ${R_{d-1}}\cap {\mathfrak{p}}$ is the only prime ideal associated to the 
 $R_{d-1}$-module $M$.
Since ${\mathfrak{p}}$ is a non-principal ideal,
by Lemma \ref{ent} the action $\alpha_{M}^{\Lambda}$ has zero entropy.
$\hfill \Box$

\section{Rigidity of equivariant maps}

 Recall that a valuation on a ring $R$ is a map 
$v : R\rightarrow {\mathbb{R}}\cup\{\infty\}$ satisfying 
the following conditions :
\begin{enumerate}
\item{$v(1) = 0$ and $v(x) = \infty$ if and only if $x = 0$,}
\item{$v(xy) = v(x) + v(y)$,}
\item{$v(x + y) \ge \mbox{Min}\{ v(x), v(y)\}.$}
\end{enumerate}

Let ${\mathfrak{p}}\subset{R_d}$ be a prime ideal. For a valuation $v$ on
${R_d}/{\mathfrak{p}}$  we define a homomorphism 
$\phi_v:{\mathbb{Z}}^d\rightarrow{\mathbb{R}}$
by $\phi_v({\bf n}) = v(u^{{\bf n}})$. The valuation $v$ is said
to be {\it discrete \/} if the image of  
$\phi_v$ is a non-trivial discrete subgroup of ${\mathbb{R}}$.
It is easy to see that $v$ is discrete if and
only if the subgroup 
$\Lambda_v = \{{\bf n} : \phi_v({\bf n}) = 0\} \subset {\mathbb{Z}}^d$
has rank $d-1$.

The following result characterizes all prime 
ideals ${\mathfrak{p}}\subset{R_d}$
with the property that  ${R_d}/{\mathfrak{p}}$ admits a discrete
valuation; 
for a proof see \cite[Theorem 2.4]{BG} and 
\cite[Corollary 1 and Theorem 8.1]{BST}. 
\begin{lemma}\label{BG}
Let $d\ge 1$, and let ${\mathfrak{p}}\subset {R_d}$ be a 
prime ideal. Then ${R_d}/{\mathfrak{p}}$ admits a discrete valuation
if and only if
it is not finitely generated as an additive group. 
\end{lemma}

{\it Definition. \/}
 Let $d\ge 1$, and let $M$ be a countable ${R_d}$-module.
A closed subgroup $H\subset X_M$ is {\it $\alpha_M$-shrinking \/}
if  there exists a ${\bf n}_{H}\in {\mathbb{Z}}^d$ such that
$$\bigcap_{j\ge 1} \alpha_{M}(j{\bf n}_{H})(H) = \{0\},$$
and $\alpha_M({\bf n})(H) = H$ whenever $<{\bf n},{\bf n}_{H}> = 0$,
where  $<.,.>$ denotes the standard inner product in ${\mathbb{R}}^d$.

\medskip
The smallest closed subgroup of $X_M$ which contains all 
$\alpha_M$-shrinking subgroups will be denoted by $X^{s}_{M}$.
It is easy to see that
if $M,N$ are ${R_d}$-modules and $\theta: X_M\rightarrow X_N$ is a 
continuous ${\mathbb{Z}}^d$-equivariant
homomorphism, then
$\theta(X_{M}^{s}) \subset X_{N}^{s}$.
\begin{prop}
\label{val}
Let $d\ge 1$, $M$ be a countable ${R_d}$-module and $F(M)$ be as defined
before. Then $X_{M}^{s} = F(M)^{\perp}$.
\end{prop}

{\it Proof.\/}
First we consider the case when $F(M) = \{0\}$.
If the above assertion is not true then
 $X^{s}_{M} = M_{1}^{\perp}$ for some 
non-zero submodule $ M_{1}$  of $M$. 
In that case we choose ${\mathfrak{p}}\in\mbox{Asc}(M_{1})$ 
and $m_{0}\in M_{1}$ with
$\mbox{ann}(m_{0}) = {\mathfrak{p}}$.
Since $F(M) = \{0\}$, ${R_d}/{\mathfrak{p}}$ is not finitely generated as an
additive group. By the previous lemma there exists a discrete valuation $v$
on  ${R_d}/{\mathfrak{p}}$. We set $$R = \{p\in {R_d} : v(p)\ge 0\},\ 
M_{0} = \{p\cdot m_{0} : v(p)\ge 1\}.$$ It is easy to see that $R$ is a
Noetherian subring and $M_{0}$ is a $R$-module.

For any $A\subset M$, let $\overline{A}$ denote the ${R_d}$-submodule
generated by $A$. We define a partial order on the set of all
$R$-submodules of $M$ by setting 
$N\le N^{'}$ if $N \subset N^{'}$ and $N^{'}\cap {\overline{N}}
= N$. It is easy to see that for any  totally ordered subset
${C}$, the $R$-module $\cup\{ N : N\in {C}\}$ is an upper bound 
for ${C}$. By Zorn's lemma there
exists a maximal element  $N_{0}$ in $\{N : N\ge M_{0}\}$.
We claim that ${\overline{N_{0}}} = M$.  Suppose this is not the case.
Let $m$ be  element of $ M - {\overline{N_{0}}}$. Since $R$ is Noetherian,
the $R$-module 
$R\cdot m \cap {\overline{N_{0}}}$ is finitely generated.
 We choose  ${\bf n}\in {{\mathbb{Z}}^d}$ with $v({\bf n}) < 0$,
and observe that
$${\overline{N_{0}}} = \bigcup_{j\ge 1} u^{j{\bf n}} \cdot N_{0}.$$
Let $B$ be a finite set which generates 
$R\cdot m \cap {\overline{N_{0}}}$ as a $R$-module.
We  find $k \ge 1$ such that 
$u^{-k{\bf n}}\cdot b\in N_{0}$ for all $b\in B$. If 
$m^{'} = u^{-k{\bf n}}\cdot m$ then
$(R\cdot m^{'} + N_{0})\cap {\overline{N_{0}}}
= R\cdot m^{'}\cap {\overline{N_{0}}} + N_{0}
= N_{0}$,
which contradicts the maximality of $N_{0}$ and proves the claim.

Now from the above claim and duality theory
we deduce that 
$$  \bigcap_{j\ge 1}\alpha(j{\bf n})(N_{0}^{\perp}) = \{ 0\},$$ 
i.e. $N_{0}^{\perp}$  is a $\alpha$-shrinking subgroup of $X_M$.
Hence $N_{0}\supset M_{1}$. As 
 $N_{0}\cap {\overline{M_{0}}}
 = M_{0}$, we deduce that 
$m_0\in M_{0}$ i.e. $p\cdot m_{0} = m_{0}$ for some
$p\in {R_d}$ with $v(p)\ge 1$. Since 
$0 = v(1)\ge {\mbox{Min}}\{v(p),v(1-p)\}$, we arrive at a  
contradiction. This completes the proof of the given assertion in the
special case considered above.

Now we consider the general case. For $a\in F(M)$, let $M_a$
denote the ${R_d}$-submodule generated by $a$, and let 
$\pi_{a} : X_{M}\rightarrow X_{M_{a}}$ denote the dual of the
inclusion map $i : M_{a}\rightarrow M$. Since $M_a$ is finitely
generated as an additive group,  $X_{M_a}$ is isomorphic with
${\mathbb{T}}^{n}\times F$, where
${\mathbb{T}}^n$ is a torus and $F$ is a finite abelian
group. Hence  $X_{M_a}$ does not admit arbitrarily small 
non-trivial subgroups. In particular, $X^{s}_{M_a} = \{0\}$.
Since $\pi_a$ is a continuous ${\mathbb{Z}}^d$-equivariant homomorphism
for every $a\in F(M)$, this shows that
$$X_{M}^{s} \subset \bigcap_{a\in F(M)}\mbox{Ker}(\pi_a)
 =  F(M)^{\perp} .$$
Let $\theta : X_{M/F(M)}\rightarrow X_{M}$ denote the dual of
the projection map from $M$ to $M/F(M)$.
Since $F(M/F(M)) = \{0\}$, 
$X_{M/F(M)}^{s} = X_{M/F(M)}$. As $\theta$ is 
a continuous ${\mathbb{Z}}^d$-equivariant homomorphism,
this implies that
$F(M)^{\perp}  = \theta(X_{M/F(M)}) 
= \theta(X_{M/F(M)}^{s})\subset X_{M}^{s}$. 
$\hfill \Box$

\bigskip
A measurable map $f : X\rightarrow Y$ between compact abelian groups
is  a {\it constant \/} if there exists a $c\in Y$ such that
$f(x) = c$ almost everywhere with respect to $\lambda_{X}$. For a 
measurable map $f : X\rightarrow Y$, by $C(f)$ we denote the
set of all $h\in X$ with the property that the map
$x\mapsto f(x+h)-f(x)$ is a constant.
We note the following elementary fact.
\begin{prop}
\label{elem}
Let  $f : X\rightarrow Y$ be a  measurable map
between compact abelian groups. Then $C(f)$ 
is a closed subgroup of $X$, and $f$
is an affine map if and only if $C(f) = X$.
\end{prop}
{\it Proof.\/}
From duality theory it follows that $f$ is an affine map if and only
if $\chi\circ f$ is an affine map for every $\chi\in {\widehat{Y}}$,
and $C(f) = \cap_{\chi} C(\chi\circ f)$.
Hence without loss of generality we may assume that $Y = {\mathbb{T}}$. It is
easy to see that $C(f)\subset X$ is a subgroup, and $C(f) = X$
whenever $f$ is an affine map. As the translation action of $X$ on
$L^2(X)$ is continuous, and the space of constant maps from $X$ to
${\mathbb{T}}$ is a closed subset of $L^2(X)$,  $C(f)$ is
closed. Since all eigenvalues of the translation action are of the
form $c\cdot \chi$ for some $c\in{\mathbb{C}}$ and $\chi\in {\widehat{X}}$,
 $C(f) = X$ only if $f$ is an affine map.
$\hfill \Box$.

\bigskip
{\it Proof of Theorem 1.1. \/}
Let $M$ and $N$ be the dual modules of $\alpha$ and $\beta$,
respectively.
For any $a\in N$, let $Y_{a}\subset Y$ denote the dual
of ${R_d}\cdot a$, and let $\beta_{a}$ denote the 
${\mathbb{Z}}^d$-action on $Y_{a}$ induced by $\beta$. If $\pi_{a}$
denotes the  factor map from $(Y,\beta)$ to
$(Y_{a},\beta_{a})$, then it is easy to see that the map
$f$ is affine if and only if $\pi_{a}\circ f$ is afiine 
for every $a\in N$. Since the dual module of each $\beta_a$ is
Noetherian, without loss of generality
we may assume that $N$ is a Noetherian ${R_d}$-module. 

We define a  map $q : X\times X\rightarrow Y$ by
$$ q(x, x^{'}) = f(x+x^{'}) - f(x^{'}).$$
For $\mu\in M(X)$, let ${\overline{\mu}}$ denote
the measure $q_{*}(\mu\times\lambda_{X})\in M(Y)$.
It is easy to see that $q$ is a ${\mathbb{Z}}^d$-equivariant map,
and ${\overline{\alpha({\bf m})(\mu)}} = 
\beta({\bf m})({\overline{\mu}})$ for all $\mu\in M(X)$ and
${\bf m}\in {{\mathbb{Z}}^d}$. We claim that 
${\overline{\mu}}\mapsto \delta_{0_Y}$ as $\mu\mapsto \delta_{0_{X}}$.
Let $\widehat{{\overline{\mu}}}$ denote the Fourier transform of 
${\overline{\mu}}$. We fix a  non-zero character $\chi\in N$, 
and note that
$${\widehat{\overline{\mu}}}(\chi) = \int \chi\ \mbox{d}{\overline{\mu}}
= \int \chi\circ q(x,x^{'})\ \mbox{d}\mu(x)\ \mbox{d}\lambda_{X}(x^{'}).$$

We set $q_{\chi} = \chi\circ q$ and $f_{\chi} = \chi\circ f$. For
$x\in X$, we denote the $L^{1}$-norm of the function
$x^{'}\mapsto q_{\chi}(x, x^{'})$ by $P(x)$. 
Since the translation action
of $X$ on $L^{1}(X)$ is strongly continuous, and
$$ | q_{\chi}(x,x^{'})| = |f_{\chi}(x + x^{'})\cdot
{\overline{f_{\chi}(x^{'})}}| \le 
 |f_{\chi}(x + x^{'})- f_{\chi}(x^{'})|,$$
it follows that $P : X\rightarrow {\mathbb{R}}$ is a continuous function
with $P(0_{X}) = 0$. Since 
$|{\widehat{\overline{\mu}}}(\chi)| \le \int P(x)\ \mbox{d}\mu(x)$, 
${\widehat{{\overline{\mu}}}}(\chi)\mapsto 0$ as 
$\mu\mapsto \delta_{0_{X}}$. This proves the claim.

We fix a $\alpha$-shrinking subgroup $H\subset X$ and define
$$\Lambda  = \{{\bf n}\in {{\mathbb{Z}}^d} :\ <{\bf n}, {\bf n}_{H}> = 0\}.$$
Since $\rm{rank}(\Lambda) = d-1$, by Lemma \ref{finite} 
the action $\beta^{\Lambda}$
has finite entropy. For $j\ge 0$, let $\mu_j$ denote the Haar measure on 
$\alpha(j{\bf n}_{H})(H)$. As $\mu_0 = \lambda_H$ is invariant under the action
$\alpha^{\Lambda}$, the measure $\overline{\mu_{0}}$ is invariant
under $\beta^{\Lambda}$. 
Since $\mu_j\mapsto \delta_{0_X}$ as $j\mapsto\infty$,
 from the above claim we deduce that
${\overline{\mu_j}}\mapsto \delta_{0_{Y}}$ as $j\mapsto\infty$.
We note that $\beta({\bf n})$ commutes
with the action $\beta^{\Lambda}$ for all ${\bf n}\in {{\mathbb{Z}}^d}$.
This implies that
 the $\Lambda$-actions $(Y, \beta^{\Lambda}, {\overline{\mu_0}})$ and 
$(Y, \beta^{\Lambda}, {\overline{\mu_j}})$ are measurably conjugate for every 
$j\ge 1$. By Lemma \ref{semia},
$$ h_{\overline{\mu_0}}(\beta^{\Lambda}) = \lim_{j\mapsto\infty}
 h_{\overline{\mu_j}}(\beta^{\Lambda}) = 0.$$

We set $X_1 = (X, \alpha^{\Lambda}, \lambda_{H}), 
X_2 = (X, \alpha^{\Lambda}, \lambda_{X})$ and
$Y_1 = (Y, \beta^{\Lambda}, {\overline{\mu_0}})$.
Since $q : X_1\times X_2\rightarrow Y_1$ is measure preserving
and $Y_1$ has zero-entropy,  the map 
$q$ is measurable with respect to the Pinsker algebra
of $X_1\times X_2$.
As Pinsker algebra of the product of two 
measure preserving ${\mathbb{Z}}^d$-actions is the product of their
Pinsker algebras
(cf. \cite[Theorem 4]{GTW}), and since $X_2$ has completely positive entropy, 
there exists a measurable map $g : X\rightarrow Y$ such that
$q(x, x^{'}) = g(x)$ almost everywhere with respect to $\lambda_{H}$.
An application of Fubini's theorem shows that for $\lambda_H$-a.e. $h\in H$,
the map $x^{'}\mapsto q(h, x^{'})$ is a constant. As 
our initial choice of $H$ was arbitrary, from the previous proposition
we deduce that $X^{s}_{M}\subset C(f)$. Hence the map
$q^1 : X\times X\rightarrow Y$ defined by
$$ q^{1}(x, x^{'}) = f(x + x^{'} ) - f(x) -f(x^{'}),$$
is invariant under the translation action of 
$X^{s}_{M}\times X^{s}_{M}$ on $X\times X$.
Let $\pi$ denote the projection map from 
$X\times X$ to $ Z = {X/X^{s}_{M}}\times  {X/X^{s}_{M}}$, and let
$q^2 : Z\rightarrow Y$ denote the  measurable map 
satisfying
$q^{1} = q^{2}\circ\pi$. 
If $\alpha_Z$ denote the ${\mathbb{Z}}^d$-action on
$Z$ induced by $\alpha$, then by 
Proposition \ref{val} the dual module of $\alpha_Z$
is isomorphic with $F(M)\times F(M)$.
Since $q^2$ is a ${\mathbb{Z}}^d$-equivariant map,
by Lemma \ref{toral} both $q^2$ and $q^{2}\circ \pi = q^1$ 
are  affine maps. As $q^1(x, x^{'}) = q^{1}(x^{'}, x)$,
this implies that for any $t\in X$,
$q^{1}(x+t, x^{'}) = q^{1}(x, x^{'} + t)$ almost everywhere
with respect to 
$\lambda_{X}\times \lambda_{X}$. Hence for any $t\in X$,
$$ f(x+t) + f(x^{'}) = f(x + t + x^{'}) - q^{1}(x+t,x^{'})
= f(x) + f(x^{'}+t),$$
from which we deduce that $X = C(f)$. By the previous proposition
 $f$ is an affine map.
$\hfill \Box$

\bigskip
An algebraic ${\mathbb{Z}}^d$-action $\alpha$ is {\it prime\/} if the 
dual module of $\alpha$ is of the form ${R_d}/{\mathfrak{p}}$ for some
prime ideal ${\mathfrak{p}}\subset {R_d}$. Prime actions can be viewed as
building blocks of algebraic ${\mathbb{Z}}^d$-actions 
(see \cite[Figure 1]{S}). 
From duality
theory it follows that $X_{{R_d}/{\mathfrak{p}}}$ is connected if 
${\mathfrak{p}}\cap {\mathbb{Z}} = \{0\}$, and it is zero-dimensional if
${\mathfrak{p}}\cap {\mathbb{Z}} \ne \{0\}$. In the zero-dimensional case
in \cite{BS} and \cite{E} it was shown that any measurable
factor map between mixing zero-entropy prime actions is an
affine map. As an application of Theorem \ref{main} we now
extend this result to all mixing zero-entropy prime actions.
  
\begin{cor}
Let $d > 1$, and  let ${\mathfrak{p}} , {\mathfrak{q}}\subset {R_d}$ be  
prime ideals such that 
the actions 
$\alpha_{{R_d}/{\mathfrak{p}}}$ and $\alpha_{{R_d}/{\mathfrak{q}}}$ are 
mixing and have 
zero-entropy. Then any measurable factor map from 
$(X_{{R_d}/{\mathfrak{p}}}, \alpha_{{R_d}/{\mathfrak{p}}})$
to $(X_{{R_d}/{\mathfrak{q}}}, \alpha_{{{R}_d}/{\mathfrak{q}}})$ 
is an affine map. 
In particular, the actions $\alpha_{{R_d}/{\mathfrak{p}}}$ and
$\alpha_{{R_d}/{\mathfrak{q}}}$ are measurably conjugate if and only if
${\mathfrak{p}} = {\mathfrak{q}}$.
\end{cor}
{\it Proof. \/}
If $d_{1}\le d$ and $\phi : {\mathbb{Z}}^{d_1}\rightarrow
{\mathbb{Z}}^{d}$ 
is an
injective homomorphism,
we define an algebraic ${\mathbb{Z}}^{d_1}$-action 
$\alpha^{\phi}$ on $X_{{R_d}/{\mathfrak{p}}}$ by setting
$\alpha^{\phi}({\bf m}) = \alpha_{{R_d}/{\mathfrak{p}}}(\phi({\bf m}))$ for all
${\bf m}\in {\mathbb{Z}}^{d_1}$. 
For any such $\phi$, let
 $\phi_{*} : {R}_{d_{1}}\rightarrow {R_d}$ denote the induced
ring homomorphism, and let 
${\mathfrak{p}}_{\phi}\subset {R}_{d_{1}}$ denote the prime ideal
$\phi_{*}^{-1}({\mathfrak{p}})$. If $M_1$ denotes the dual module 
of $\alpha^{\phi}$
then it is easy to see that $\mbox{Asc}(M_1) = \{{\mathfrak{p}}_{\phi}\}$.
By Lemma \ref{ent} the action $\alpha^{\phi}$ has zero-entropy if
${\mathfrak{p}}_{\phi}$ is non-principal, and it has completely positive 
entropy if ${\mathfrak{p}}_{\phi}$ is principal.
Let $d_{0}$ be the smallest integer such that for some injective
homomorphism $\phi_{0} : {\mathbb{Z}}^{d_0}\rightarrow
{\mathbb{Z}}^{d}$ 
the ideal
 ${\mathfrak{p}}_{\phi_0}\subset {R}_{d_0}$ is non-principal.
 We set $\Lambda = \phi({\mathbb{Z}}^{d_0})$, and note that the action
$\alpha_{{R_d}/{\mathfrak{p}}}^{\Lambda}$ has zero-entropy but all 
lower rank subactions
of $\alpha_{{R_d}/{\mathfrak{p}}}^{\Lambda}$ have completely positive 
entropy. Since 
$(X_{{R_d}/{\mathfrak{q}}}, \alpha_{{{R}_d}/{\mathfrak{q}}})$ 
is a factor of 
$(X_{{R_d}/{\mathfrak{p}}},\alpha_{{R_d}/{\mathfrak{p}}})$, the same
is 
true for the
action $\alpha_{{R_d}/{\mathfrak{q}}}^{\Lambda}$. By Theorem 1.1, 
$f$ is an affine map.
Hence if $f$ is a measurable conjugacy then the actions 
$\alpha_{{R_d}/{\mathfrak{p}}}$ and $\alpha_{{R_d}/{\mathfrak{q}}}$
are algebraically conjugate. By duality, this happens only if
${R_d}/{\mathfrak{p}}$ and ${R_d}/{\mathfrak{q}}$ 
are isomorphic ${R_d}$-modules,
 i.e., ${\mathfrak{p}} = {\mathfrak{q}}$.

$\hfill \Box$




\end{document}